\documentclass[a4paper,11pt]{amsart}

\usepackage{enumerate}

\usepackage{upref}

\usepackage[colorlinks,backref]{hyperref}                   

\usepackage{amssymb,amsmath}

\newtheorem{theorem}{Theorem}[section]
\newtheorem{lemma}[theorem]{Lemma}

\newtheorem{definition}{Definition}[section]

\theoremstyle{remark} 
\theoremstyle{definition} 

\numberwithin{equation}{section}

\newif\ifcomment \commentfalse \def\commentON{\commenttrue}
 \long\outer\def\BC#1\EC{\ifcomment \sloppy \par
\# \ldots\dotfill {\em #1} \dotfill \# \par \fi } \commentON

\newcommand{\remove}[1]{}

\newcommand{\eps}{\ensuremath{\varepsilon}}

\newcommand{\R}{{{\mathbb{R}}}}

\newcommand{\Div}{\mathrm{div}}

\newcommand{\dx}{\, dx}

\newcommand{\dy}{\, dy}

{%

\begin{enumerate}}%
{\end{enumerate}
}

\begin{document}

\title[Dissipation for a non-convex gradient flow problem]{Dissipation for a non-convex gradient flow problem
of a Patlack-Keller-Segel type for densities on $\R^d$, $d \geq 3$}

\author[{E. A. Carlen}]{Eric A. Carlen}
\address[ Eric A. Carlen]{\newline
	Department of Mathematics\newline
        Rutgers University \newline
	Piscataway, NJ, USA}
\email{carlen@math.rutgers.edu}
\urladdr{http://www.math.rutgers.edu/~carlen/}

\author[{S. Ulusoy}]{Suleyman Ulusoy}
\address[ Suleyman Ulusoy]{\newline
	Department of Mathematics and Natural Sciences\newline
        American University of Ras Al Khaimah\newline
	PO Box 10021, Ras Al Khaimah, UAE}
\email{suleyman.ulusoy@aurak.ac.ae}
\urladdr{https://aurak.ac.ae/en/dr-suleyman-ulusoy/}




\date{\today}

\begin{abstract} We study an evolution equation that is the gradient flow in the $2$-Wasserstien metric of a non-convex functional for densities in $\R^d$ with $d\geq 3$. Like the Patlack-Keller-Segel system on $\R^2$, this evolution equation features a competition between the dispersive effects of diffusion, and the accretive effects of a concentrating drift. We determine a parameter range in which the diffusion dominates, and all mass leaves any fixed compact subset of $\R^d$ at an explicit polynomial rate.

\end{abstract}

\maketitle


\section{Introduction}
\label{sec:intro}

\subsection{The gradient flow equation:}

Lieb's  sharp form of Hardy-Littlewood-Sobolev(HLS) inequality \cite{24} states that for a nonnegative measurable function $f$ on $\R^d,$ and all $0< \lambda < d,$
\begin{equation}\label{HLS}
\frac{\int_{\R^d} \int_{\R^d} \frac{f(x)f(y)}{|x-y|^{\lambda}} \dx \dy }{||f||_p^2} \leq \frac{\int_{\R^d} \int_{\R^d} \frac{h(x)h(y)}{|x-y|^{\lambda}} \dx \dy }{||h||_p^2} =: C_{HLS}
\end{equation}
where $h = (1+|x|^2)^{(\lambda -2d)/2}$ and $p = \frac{2d}{2d-\lambda}.$

We focus on the case $\lambda = d-2$ and for $\alpha > 0$ define the functional
\begin{equation}\label{enfnc}
E_\alpha[f] := ||f||_{\frac{2d}{d+2}}^2 - \frac{\alpha}{C_{HLS} }\int_{\R^d} f(x) \left[  (-\Delta)^{-1} f \right](x) \dx .
\end{equation}
Note that In for all $\alpha > 0$, $E_\alpha$ is {\em difference} of two convex functions. Moreover,
 for $\alpha\in (0,1)$, $E_\alpha$ is strictly positive,  though for $f_{c,s,x_0}(x) := h(x/s - x_0)$, $s>0$, $x_0\in \R^d$, an HLS optimizer,
$$E_\alpha(f_{c,s,x_0}) =  s^d(1-\alpha)c^2 ||h||_{\frac{2d}{d+2}}^2\ .$$
Hence for $\alpha <1$, and $s\downarrow 0$, $E_\alpha(f_{c,s,x_0}) \downarrow 0$.  It follows that  $\alpha  < 1$, the functional $E_\alpha$ has no non-zero minimizers.
By the same computation, for $\alpha > 1$, $E_\alpha$ is not convex and not bounded below, and as $s\uparrow \infty$, $E_\alpha(f_{c,s,x_0}) \downarrow -\infty$.  In fact, for all $\alpha > 0$, $E_\alpha$ is {\em difference} of two convex functions

Note that $E_1[f] \geq 0$ by the HLS inequality. For $\alpha > 1$, $E_\alpha[f] $ is not bounded below, while for
 $\alpha < 1$ it is bounded from below. Then we have
 \begin{eqnarray}\label{IC2}
 E_\alpha[f]   &=& E_1[f]  +  \frac{(1-\alpha)}{C_{HLS} }\int_{\R^d} f(x) \left[  (-\Delta)^{-1} f \right](x) \dx .\nonumber\\
 &=&  \alpha E_1[f]  + (1-\alpha) ||f||_{\frac{2d}{d+2}}^2 \ .
 \end{eqnarray}

This functional arises in one of several natural ways to generalize the Patlack-Keller-Segel system to more than two spatial dimensions. The evolution equation we consider is the gradient flow for $E_\alpha$
in the two Wasserstein metric. This gradient flow equations can be written as
\begin{equation}\label{grflow}
\frac{\partial f}{\partial t} = \Div \left(f \nabla \left(\frac{\delta E_\alpha}{\delta f}\right) \right),
\end{equation}
where $\frac{\delta E_\alpha}{\delta f}$ is the first variation of the energy with respect to the $L^2-$metric.

In writing this out more explicitly, it will be convenient to define $\kappa = \alpha/C_{HLS}$.  We obtain:
\begin{equation}\label{Meqn}
 \frac{\partial f}{\partial t} = \left(\frac{d-2}{d}\right) ||f||_{\frac{2d}{d+2}}^{\frac{4}{d+2}} \Big\{  \Delta\left(f^{\frac{2d}{d+2}}\right)    - \Div \left(  f \nabla c  \right)  \Big\},
\end{equation}
where $\nabla c$ is given by
\begin{equation}\label{nabc}
\nabla c = \frac{d \kappa } {(d-2)||f||_{\frac{2d}{d+2}}^{\frac{4}{d+2}}}  \nabla \Big(   \Big[  (-\Delta)^{-1}  f  \Big](x) \Big).
\end{equation}
We consider  equation (\ref{Meqn}) for initial data
\begin{equation}\label{ic}
f(t=0, x) = f_0(x),
\end{equation}
that satisfies
\begin{equation}\label{IC}
0 \leq f_0(x), \quad (1 + |x|^2) f_0(x) \in L^1(\R^d), \quad E_\alpha[f_0] < \infty.
\end{equation}

For more background on this equation and information on the relation to the Patlack-Keller-Segel system,
see the paper  \cite{ulusoy}, where a number of results on existence and blow-up are proved. The non-convexity of the functional $E_\alpha$ is the source of interesting features in the study of this gradient flow problem.

Before stating our result, we recall some relevant facts about the original two dimensional problem. In the theory of this problem as developed by Dolbeault and Perthame \cite{17},  the logarithmic Hardy-Littlewood-Sobolev  (log-HLS) inequality plays the role corresponding to (\ref{enfnc}). The sharp log-HLS inequality \cite{B, CL} states that for all  $f\geq 0$ on $\R^2$ such that $f \log f$ is integrable,
\begin{equation}\label{lhls}
\int_{\R^2}f \log f{\rm d}x  + \frac{2}{M}\int_{\R^2\times \R^2} f(x) \log |x-y| f(y) {\rm d}x{\rm d}y \geq M(1 +\ log \pi -\log M)
\end{equation}
where
$M = \int_{\R^2} f {\rm d}x$ is the {\em total mass}. Thus, the log HLS inequality differs form the HLS
inequality in that the sharp constant in it depends on the mass. Building on earlier work of J\"ager and Luckhaus \cite{JL},
Dolbeault and Perthame showed that the mass $M = 8\pi$ is critical in two dimensions. For initial data $f_0$ with
$ \int_{\R^2} f_0 {\rm d}x < 8\pi$, solutions $f(t,x)$ exist for all time, and $\int f(t,x) \log f(t,x){\rm d}x$ is
uniformly bounded in $t$. J\"ager and Luckhaus proved existence of solutions and obtained this global existence
result for solutions with $ \int_{\R^2} f_0 {\rm d}x < M_0$, for an explicit $M_0$ smaller than $8\pi$. For initial data
with $f_0$ with  $M := \int_{\R^2} f_0 {\rm d}x > 8\pi$, Dolbeault and Perthame showed
solutions with finite second moment  blow up in a finite time depending only on $M - 8\pi$ and the
second moment of the initial data. The mass $M = 8\pi$  makes the left side of (\ref{lhls}) exactly the functional for which the $2$-dimensional  Patlack-Keller-Segel system is gradient flow.

In the higher dimensional case, the constant $\alpha$ plays the role of the mass parameter. In this case,  for
 $\alpha =1$,  the driving functional for our equation is precisely the functional that is non-negative by the sharp HLS inequality. Therefore, and one may conjecture that $\alpha =1$ is the critical parameter value for our equation.
This has indeed been proved in previous work of one of the authors:
When $\alpha < 1$,
solutions for a natural class of initial data exist for all time, but for $\alpha>1$ this need not be the case.

The next question is whether for $\alpha < 1$, the solutions have a particular, simple asymptotic behavior. In the $2$-dimensional case, this problem was investigated in \cite{BCP}, which provides an affirmative answer.
Theorem~\ref{decaythm} below provides such a result for our higher dimensional equation, but only for small $\alpha$, as in \cite{BDEF} for related equation in two dimensions.
 It is an open problem  in our case to remove the small $\alpha$ restriction, and to prove  that the solutions approach self-similar scaling solutions with a particular profile.  However, we do prove that for small $\alpha$, the evolution is diffusion dominated, with all mass leaving any fixed compact set at an explicitly computable rate. With this background explained, we state our main result:

 \begin{theorem}\label{decaythm} There is an  explicitly computable $0 < \alpha_0 < 1$ such that for all
 $\alpha \in (0,\alpha_0)$, a solution of our equation with initial data satisfying
 (\ref{IC2})  satisfies
 $$E_\alpha[f(\cdot,t)] \leq Ct^{-\frac{d-2}{d}}\ $$
 where $C$ is a constant depending only on the initial mass $\int_{\R^d}f_0(x) \, {\rm d}x$
 and free energy $E_\alpha(f_0)$.
 \end{theorem}

 Note that by (\ref{IC2}), this means that for such solutions $f(x,t)$,
 $$ ||f(\cdot, t) ||_{\frac{2d}{d+2}}^2  \leq \frac{C}{1-\alpha_0} t^{-\frac{d-2}{d}}\ .$$
 Since the total mass is conserved, this means that such solutions spreads out, with the mass in any fixed
  compact subset of $\R^d$ decaying to zero at a polynomial rate. Thus, for
  $\alpha < \alpha_0$, the diffusion dominates the concentrating drift.

The key to proving this is the following lemma:

\begin{lemma}\label{Lqlem}

\begin{equation}\label{g4}
\frac{d}{dt}\int f^q(x, t) \, dx \leq (q-1)\left[-\frac{8(d-2)q}{(d+1)(q+\frac{d-2}{d+2})^2} \frac{1}{C_{GNS}^{\frac{2(q+1)(d+2)}{q(d+2)+d-2}}}  + \kappa \right] ||f||_{q+1}^{q+1}.
\end{equation}
where $C_{GNS}$ is the sharp constant in the Gagliardo-Nirenberg-Sobolev inequality
\begin{equation}\label{g3}
||u||_{\frac{4d}{q(d+2)+d-2}}^{\frac{4}{(q+1)(d+2)}} ||\nabla u||_2^{\frac{q(d+2)+d-2}{(q+1)(d+2)}} \geq \frac{1}{C_{GNS}} ||u||_{\frac{2(q+1)(d+2)}{q(d+2)+d-2}}.
\end{equation}
\end{lemma}

We apply the lemma for $q = 2d/(d-2)$, so that $q+1 = (3d+2)/(d+2)$. Then  the right hand side in (\ref{g4}) becomes
$$
\frac{d-2}{d+2}\left[-\frac{16d(d^2-4)}{(d+1)(3d-2)^2} \frac{1}{C_{GNS}^{2(3d+2)/(3d-2)} }
+ \kappa \right]   =: -K
$$
Then for
$$\alpha_0 =  \frac{16d(d^2-4)}{(d+1)(3d-2)^2} \frac{C_{HLS}}{C_{GNS}^{2(3d+2)/(3d-2)} }\ ,$$
$\int_{\R^d }f^{2d/(d+2)}(x, t) \, dx$ is monotone decreasing for all $\alpha <  \alpha_0$.  Integrating, we have that for any $t_0 < t_1$,
$$ \int_{\R^d } f^{\frac{2d}{d+2}} (x, t_1) \, dx  - \int_{\R^d } f^{\frac{2d}{d+2}} (x, t_0)  \, dx  \leq  -K \int_{t_0}^{t_1} \int_{\R^d } f^{\frac{3d+2}{d+2}}(x,t) \, dx \, dt $$
Therefore,
$$\frac{1}{t_1- t_0}  \int_{t_0}^{t_1}\int_{\R^d } f^{\frac{3d+2}{d+2}} (x,t) \, dx \,dt  \leq \frac{1}{(t_1-t_0)K} \int_{\R^d }f^{\frac{2d}{d+2}} (x, t_0)  \, dx \ .$$
By (\ref{IC2}),
$$\int_{\R^d } f^{\frac{2d}{d+2}} (x, t_0) \, dx  \leq \left(\frac{E_\alpha[f(\cdot,t_0)]}{1-\alpha}\right)^{\frac{d}{d+2}}\ .$$
Then since a minimum cannot exceed an average, there is some $t \in [t_0,t_1]$ with
\begin{eqnarray}
\int_{\R^d }f^{\frac{3d+2}{d+2}}(x,t) \, dx  &\leq & \frac{1}{t_1- t_0}\int_{t_0}^{t_1}  \int_{\R^d }f^{\frac{3d+2}{d+2}}(x,t) \, dx  \, dt \nonumber\\
&\leq&  \frac{1}{(t_1-t_0)K}   \left(\frac{E_\alpha[f(\cdot,t_0)]}{1-\alpha}\right)^{\frac{d}{d+2}}  \ .
\end{eqnarray}
Now define $\beta$ so that
$$\frac{2d}{d+2}   = \beta  \frac{3d+2}{d+2}   + (1-\beta) 1\ ,$$
which means that $\beta = \frac{d-2}{2d}$.  Then by H\"older's inequality,
$$\int_{\R^d} f^{\frac{2d}{d+2}} \, dx  \leq \left( \int_{\R^d} f^{\frac{3d+2}{d+2}} \, dx \right)^\beta \left( \int_{\R^d} f \, dx \right)^{1-\beta}\ .$$
Letting $M$ denote the conserved mass $M = \int_{\R^d} f \,dx$, we now have that for some $t\in [t_0,t_1]$,
\begin{eqnarray}
\int_{\R^d} f^{\frac{2d}{d+2}}(x,t) \, dx  &\leq& \left(\int_{\R^d} f^{\frac{3d+2}{d+2}}(x,t) dx \right)^{\frac{d-2}{2d}} M^{\frac{d+2}{2d}} \nonumber\\
&\leq& \left(\frac{1}{(t_1-t_0)K}   \left(\frac{E_\alpha[f(\cdot,t_0)]}{1-\alpha}\right)^{\frac{d}{d+2}} \right)^{\frac{d-2}{2d}} M^{\frac{d+2}{2d}}
\end{eqnarray}

Then, using the monotonicity of the free energy,  $ E_\alpha[f(\cdot,t_1)]^{\frac{d}{d+2}} \leq  E_\alpha[f(\cdot,t)]^{\frac{d}{d+2}}  \leq  \int_{\R^d} f^{\frac{2d}{d+2}}(x,t) \,dx$. Thus we have
\begin{equation}\label{iter}
E_\alpha[f(\cdot,t_1)]  \leq \left(\frac{1}{(t_1-t_0)K} \right)^{\frac{d^2-4}{2d^2}}  M^{1/2}\left(\frac{E_\alpha[f(\cdot,t_0)]}{1-\alpha} \right)^{\frac{d-2}{2d}}\ .
\end{equation}
Now for arbitrary $t>0$, we may choose $t_1 = t$ and $t_0 = t/2$. Again using  the monotonicity of the free energy, we obtain
$$E_\alpha[f(\cdot,t)]  \leq \left(\frac{2}{K t} \right)^{\frac{d^2-4}{2d^2}}  M^{1/2}\left(\frac{E_\alpha[f_0]}{1-\alpha} \right)^{\frac{d-2}{2d}}\ .$$
Now we can use this estimate in (\ref{iter}) once more with $t_1 = t/2$ and $t_0 = t/4$ to improve the bound.  This leads to the decay rate claimed in the theorem.

\section{Proof of the main lemma}

\begin{equation}\label{g1}
\begin{split}
\frac{d}{dt}\int f^q(x, t) \, dx &= q \int f^{q-1}(x, t)f_t(x, t) \, dx \\
				       &= Aq\int f^{q-1}(x, t) [\Delta(f^{\frac{2d}{d+2}}) - \Div (f \nabla c) ](x, t) \, dx,
\end{split}
\end{equation}
where $A := (\frac{d-2}{d})||f||_{\frac{2d}{d+2}}^{\frac{4}{d+2}}$, and $\Delta c = - \frac{\kappa}{A} f$. A straightforward calculation yields that
\begin{equation}\label{g2}
\frac{d}{dt}\int f^q \, dx =  I  + II,
\end{equation}
where
$$ I = -\frac{8(d-2)q(q-1)}{(d+1)(q+\frac{d-2}{d+2})^2}||f||_{\frac{2d}{d+2}}^{\frac{4}{d+2}} \int |\nabla(f^{\frac{q+\frac{d-2}{d+2}}{2}})|^2 \, dx,$$
$$ II = (q-1)\kappa \int f^{q+1} \, dx.$$
Let $u=f^{\frac{q+\frac{d-2}{d+2}}{2}} \iff f = u^{\frac{2}{q+\frac{d-2}{d+2}}}$. Then we rewrite $I$ and $II$ in terms of $u$ as follows
$$ I =  -\frac{8(d-2)q(q-1)}{(d+1)(q+\frac{d-2}{d+2})^2} ||u||_{\frac{4d}{q(d+2)+d-2}}^{\frac{8}{q(d+2)+d-2}} ||\nabla u||_2^2,$$
$$ II = (q-1)\kappa ||u||_{\frac{2(q+1)(d+2)}{q(d+2)+d-2}}^{\frac{2(q+1)(d+2)}{q(d+2)+d-2}}.$$

We now write down the corresponding Gagliardo-Nirenberg-Sobolev type inequality. For the correct powers to be chosen one can take the $\frac{q(d+2)+d-2}{2(q+1)(d+2)}$ power of the norms in $I$ and $II$.
\begin{equation}\label{g3}
||u||_{\frac{4d}{q(d+2)+d-2}}^{\frac{4}{(q+1)(d+2)}} ||\nabla u||_2^{\frac{q(d+2)+d-2}{(q+1)(d+2)}} \geq \frac{1}{C_{GNS}} ||u||_{\frac{2(q+1)(d+2)}{q(d+2)+d-2}}.
\end{equation}
Using \eqref{g3} in \eqref{g2} yields the result.

%
%

\section{Appendix}

Here we provide the details of the existence of weak solutions with sufficient regularity to justify our calculations.
We  closely follow Bain and Lui \cite{BL}, and for this reason only sketch some arguments and make specific references to their paper for details.

Our energy functional can be rewritten in the form

\begin{equation}\label{a1}
E_{\alpha}[f] = ||f|||_{\frac{2d}{d+2}}^2 - \kappa \frac{c_d}{2}\iint_{\R^d \times \R^d} \frac{f(x)f(y)}{|x-y|^{d-2}} \, dy \, dx.
\end{equation}
Here we recall that $\kappa = \frac{\alpha}{C_{HLS}}$ and we also note that
\begin{equation}\label{a2}
c(x, t) = c_d \int_{\R^d} \frac{f(y, t)}{|x-y|^{d-2}} \, dy,
\end{equation}
where
\begin{equation}\label{a3}
c_d = \frac{1}{d(d-2)\alpha_d}, \quad \alpha_d = \frac{\pi^{d/2}}{\Gamma(d/2+1)},
\end{equation}
with $\alpha_d$ being the volume of the $d-$dimensional unit ball. Following Lemma 2.7 of    \cite{BL} we can prove that

\begin{lemma}\label{lemma2.7}
If $f \in L_{+}^1 \cap  L^{\frac{3d}{d+2}}(\R^d)$ and $c$ is given by \eqref{a2}. Then
\begin{equation}
||f |\nabla c|^2 ||_1   \leq   C ||f||_{\frac{3d}{d+2}}^3 < \infty.
\end{equation}
\end{lemma}
We now give the definition of a \emph{weak solution} and a \emph{weak entropy solution} of \eqref{Meqn}-\eqref{nabc} with initial data $f_0$ satisfying \eqref{ic} and \eqref{IC}.

\begin{definition}\label{def2.4}
Let $f_0$ be the initial data satisfying \eqref{IC} and $T \in (0, \infty)$. Let $c$ be the chemical concentration associated with $f$ and given by \eqref{a2}. $f$ is called a weak solution to   \eqref{Meqn}-\eqref{nabc} with initial data $f_0$ satisfying \eqref{ic} and \eqref{IC} if

\begin{enumerate}[(i)]

\item Regularity:

\begin{equation}\label{a4}
f \in L^2(0, T; L_{+}^1 \cap L^{\frac{2d}{d+2}}(\R^d)),
\end{equation}

\begin{equation}\label{a5}
 \partial_t f   \in L^{p_2} (0, T; W_{loc}^{-1, p_1}(\R^d)), \quad \text{for some $p_1, p_2 \geq 1$}
\end{equation}

\item  For all $\psi \in C_0^\infty(\R^d)$ and any $0<t<\infty$,
\begin{equation}\label{a6}
\begin{split}
&\int_{\R^d} \psi f(\cdot, t) \, dx - \int_{\R^d} \psi f_0 \, dx   =  \left(   \frac{d-2}{d}\right) ||f||_{\frac{2d}{d+2}}^{4/{d+2}} \int_0^t \int_{\R^d} \Delta \psi f^{\frac{2d}{d+2}} \, dx \, ds \\
&  -  \kappa \frac{c_d(d-2)}{2} \int_0^t \iint_{\R^d \times \R^d} \frac{\left[ \nabla \psi(x) - \psi(y) \right] \cdot (x-y)}{|x-y|^2}  \frac{f(x, s) f(y, s)}{|x-y|^{d-2}} \, dy \, dx \, ds
\end{split}
\end{equation}

\end{enumerate}

\end{definition}

\begin{definition}\label{def2.6}
The weak solution is also a \emph{weak entropy solution} if $f$ satisfies  additional regularity properties
\begin{equation}\label{a7}
\nabla f^{\frac{3d-2}{2d+4}}  \in L^2(0, T; L^2(\R^d))
\end{equation}
\begin{equation}\label{a8}
u \in L^3(o, T; L^{\frac{3d}{d+2}}(\R^d))
\end{equation}
and it  satisfies  the following energy dissipation inequality
\begin{equation}\label{a9}
E_{\alpha}[f(\cdot, t)] + \int_0^t \int_{\R^d }  f \left|   \nabla \left( 2||f||_{\frac{2d}{d+2}}^{4/{d+2}} f^{\frac{d-2}{d+2}}    - c \right) \right|^2 \, dx    \, ds \leq  E_{\alpha}[f_0(\cdot)],
\end{equation}
for any  $t>0$.

\end{definition}

\begin{theorem}\label{thm2.11}

Let $d \geq 3$. If $\alpha < 1$ there exists a weak entropy solution to  \eqref{Meqn}-\eqref{nabc} on $[0, \infty)$  with  initial data $f_0$ satisfying \eqref{IC} satisfying the energy inequality \eqref{a9}. If $\alpha > 1$  then there are initial data satisfying the \eqref{IC} with a negative free energy $E[f_0]$. Moreover, if  $f_0$ is such an initial data and $f$ is a weak entropy solution on $[0, T_{max}$ with initial condition $f_0$, then $T_{max} < \infty$  and $||f||_{\frac{2d}{d+2}} \to \infty$  as $t \nearrow T_{max}$.

\end{theorem}

Theroem~\ref{thm2.11} reinforces our earlier assertion that $\alpha = 1$ is the critical  parameter for the underlying system considered in this paper. As usual, the approach is to regularize the equation \eqref{Meqn} and then passing to the limit as $\eps \to 0$, we conclude the existence of the  weak entropy solution satisfying the energy inequality \eqref{a9}.  Part of this theorem was proved in \cite{ulusoy}. For other parts  of the estimates and the results we follow the ideas in \cite{BL}.

Our regularized equation is

\begin{equation}\label{aa1}
\left\{
\begin{split}
&\frac{\partial f_{\eps}}{\partial t} =  \left( \frac{d-2}{d}\right)||f_{\eps}||_{\frac{2d}{d+2}}^{\frac{4}{d+2}} \Delta\left( f_{\eps}^{\frac{2d}{d+2}}\right)\\
& \qquad \qquad \qquad  - \kappa \,  \Div \left( f_{\eps} \nabla c_{\eps}\right) + \eps \Delta f_{\eps}, \, x \in \R^d, t > 0, \\
&\Delta c_{\eps} = - J_{\eps} \ast f_{\eps}, \, x \in \R^d, t \geq 0, \\
&f_{\eps}(x, 0) = f_{0\eps}, \, x \in \R^d.
\end{split}\right.
\end{equation}
In \eqref{aa1}, $J_{\eps} := \frac{1}{\eps^d} J(\frac{x}{\eps}),$ with $J(x) = \frac{1}{\alpha_d}(1+|x|^2)^{-(d+2)/2}$  and $\int_{\R^d} j_{\eps}(x) \, dx = 1.$ This implies that
\begin{equation}\label{aa2}
c_{\eps}(x,t) = c_d \int_{\R^d} \frac{f_{\eps}(y, t)}{\left( |x-y|^2 + \eps^2\right)^{(d-2)/2}} \, dy,
\end{equation}
where
\begin{equation}\label{aa3}
c_d = \frac{1}{d(d-2)\alpha_d}, \quad \alpha_d = \frac{\pi^{d/2}}{\Gamma(\frac{d}{2}+1)}.
\end{equation}
$\alpha_d$ is the volume of the $d-$dimensional unit ball. Here $f_{0\eps}\in C^{\infty}(\R^d)$ is a sequence of approximation for $f_{0}(x)$ and can be constructed and satisfies that $\exists \delta > 0$ such that for all $0 < \eps < \delta,$
\begin{equation}\label{aa4}
\begin{split}
f_{0\eps} &> 0, \\
f_{0\eps} & \in L^r(\R^d), \, \forall r \geq 1,\\
||f_{0\eps}(x, 0)||_1 &= ||f_0(x)||_1\\
\int_{\R^d} |x|^2 f_{0\eps}(x) \, dx &\to \int_{\R^d} |x|^2 f_{0}(x) \, dx \quad \text{as} \quad \eps \to 0.
\end{split}
\end{equation}
On the other hand, if $f_0 \in L^q$  for some $q,$ then
\begin{equation}\label{aa5}
f_{0\eps}  \to f_0 \quad \text{in} \quad L^q \quad \text{as} \quad \eps \to 0.
\end{equation}
From parabolic theory, for any fixed $\eps > 0$, \eqref{aa1} has a global smooth positive solution $f_{\eps}$ with the regularity $\forall r \geq 1,$
\begin{equation}\label{aa6}
f_{\eps} \in L^{\infty}(0, T; L^r(\R^d)) \cap L^{r+1}(0, T; L^{r+1}(\R^d)).
\end{equation}
Below $\psi_R(x)$ is defined as follows: Take a cut-off function $0 \leq \psi_1(x) \leq 1$ defined as
\begin{equation}\label{aa7}
\psi_1(x) = \begin{cases}
1, & \mbox{if } |x| \leq 1 \\
0, & \mbox{if } |x| > 2,
\end{cases}
\end{equation}
where $\psi_1(x) \in C_0^{\infty}(\R^d).$ Define $\psi_R(x) := \psi_1(\frac{x}{R}),$ as $R \to \infty, \, \psi_R  \to 1 $  then there exists constants $C_1, C_2$ such that
\begin{equation}\label{aa8}
|\nabla \psi_R(x)| \leq \frac{C_1}{R}, \quad |\Delta \psi_R(x)| \leq \frac{C_2}{R^2}, \quad \text{for} \quad x \in \R^d.
\end{equation}
For future reference let us define
\begin{equation}\label{aa9}
C(d, f_\eps) := \left( \frac{d-2}{d}\right) ||f_\eps||_{\frac{2d}{d+2}}^{4/{d+2}}.
\end{equation}
For $q \geq \frac{2d}{d+2}$ we multiply \eqref{aa1} by $q f_{\eps}^{q-1} \psi_R(x)$, which yields
\begin{equation}\label{aa10}
\begin{split}
&\int f_\eps^q(x, t) \psi_R(x) \, dx - \int f_{0\eps}^q(x) \psi_R(x) \, dx \\
&\quad + \frac{C(d, f_\eps)q(q-1)(\frac{2d}{d+2})}{\left( \frac{\frac{2d}{d+2}+ q -1}{2} \right)^2} \int_0^t \int \left| \nabla f_{\eps}^{\frac{\frac{2d}{d+2}+ q -1}{2}}   \right|^2 \psi_R(x) \, dx \, ds \\
&\quad + \eps \frac{4(q-1)}{q} \int_0^t \int  \left| \nabla f_\eps^{q/2}\right|^2 \psi_R(x) \, dx \, ds \\
& \quad = (q-1)\kappa \int_0^t \int f_\eps^q (J_\eps \ast f_\eps) \psi_R(x) \, dx \, ds \\
&\qquad + \frac{C(d, f_\eps) q (\frac{2d}{d+2})}{\left(\frac{2d}{d+2} + q - 1 \right)}  \int_0^t \int f_\eps^{\frac{2d}{d+2}+q-1} \Delta \psi_R(x) \, dx \, ds \\
& \qquad + \kappa \int_0^t \int f_\eps^q \nabla c_\eps \cdot \nabla \psi_R(x) \, dx \, ds.
\end{split}
\end{equation}
where the space integration is done over $\R^d.$ First we have the following estimate:

\begin{equation}\label{aa11}
\begin{split}
\int f_\eps^q \nabla c_\eps \cdot \nabla \psi_R(x) \, dx & \leq \frac{C}{R} \int f_\eps^q |\nabla c_\eps| \, dx \\
                                                         & \leq \frac{C}{R} ||f_\eps^{\frac{2d}{d+2}}||_{r_1} ||\nabla c_\eps ||_{r_2}\\
                                                         &\leq   \frac{C}{R} ||f_\eps^{\frac{2d}{d+2}}||_{r_1} ||f_\eps||_{r_3}||\frac{x}{|x|^d}||_{L_{w}^{d/{d-1}}}\\
                                                         &\leq   \frac{C}{R} ||f_\eps||_{\frac{d(q+1)}{d+1}}^{q+1},
\end{split}
\end{equation}
where the exponents satisfy $\frac{1}{r_1}+\frac{1}{r_2} = 1, \quad \frac{1}{r_3}+ \frac{d-1}{d} = 1 + \frac{1}{r_2}$ and $qr_2 = r_3.$  Hence, we obtain

\begin{equation}\label{aa12}
\kappa \int_0^t \int f_\eps^q \nabla c_\eps \cdot \nabla \psi_R(x) \, dx \, ds \leq \frac{C(||f_{0\eps}||_1, t, \kappa)}{R} \int_0^t ||f_\eps||_{q+1}^{q+1} \, ds.
\end{equation}

We also have

\begin{equation}\label{aa13}
\frac{C(d, f_\eps)q(\frac{2d}{d+2})}{(\frac{2d}{d+2}+ q-1)} \int_0^t \int f_{\eps}^{\frac{2d}{d+2}+q-1} \Delta \psi_R \, dx \, ds  \leq \frac{C}{R^2} \int_0^t ||f_\eps||_{\frac{2d}{d+2}+q-1}^{\frac{2d}{d+2}+q-1} \, ds,
\end{equation}
where H\" older inequality with $\frac{2d}{d+2}+q-1 \geq \frac{4d}{d+2}-1 > 1$ is used. By using \eqref{aa6} and the dominated convergence theorem after taking the limit as $R \to \infty$ in \eqref{aa10} we get
\begin{equation}\label{aa14}
\begin{split}
&\int f_\eps(t)^q \, dx - \int f_{0\eps}^q(x) \, dx \\
& \quad +  \frac{C(d, f_\eps)q(q-1)(\frac{2d}{d+2})}{\left( \frac{\frac{2d}{d+2}+ q -1}{2} \right)^2} \int_0^t\int \left| \nabla f_\eps^{\frac{\frac{2d}{d+2}+q-1}{2}}\right|^2 \, dx \, ds \\
&\quad + \eps \frac{4(q-1)}{q}\int_0^t \int \left| \nabla f_\eps^{q/2} \right|^2 \, dx \, ds  \\
& \qquad = (q-1)\kappa \int_0^t \int f_\eps^q (J_\eps \ast f_\eps) \, dx \, ds.
\end{split}
\end{equation}
Now, we take the time derivative of the above equation, for any $t > 0,$
\begin{equation}\label{aa15}
\begin{split}
&\frac{d}{dt}\int f_\eps^q + \frac{C(d, f_\eps)q(q-1)(\frac{2d}{d+2})}{\left( \frac{\frac{2d}{d+2}+ q -1}{2} \right)^2} \int \left| \nabla f_\eps^{\frac{\frac{2d}{d+2}+q-1}{2}}\right|^2 \, dx \\
& + \eps \frac{4(q-1)}{q} \int  \left| \nabla f_\eps^{q/2} \right|^2 \, dx = (q-1)\kappa  \int f_\eps^q (J_\eps \ast f_\eps) \, dx.
\end{split}
\end{equation}
We can easily follow the estimates in Steps 1-6 of \cite{BL} modifies to  our case. The following estimate will be useful.
\begin{equation}\label{aa16}
\begin{split}
\int f_\eps^q (J_\eps \ast f_\eps) \, dx &\leq ||f_\eps^q||_{\frac{q+1}{q}} ||J_\eps \ast f_\eps||_{q+1} \leq ||f_\eps||_{q+1}^q ||f_\eps||_{q+1}\\
& S_d^{-1} ||\nabla f_\eps^{(\frac{2d}{d+2}+q-1)/2}||_2^2 ||f_\eps||_{\frac{2d}{d+2}}^{2-\frac{2d}{d+2}}.
\end{split}
\end{equation}
For initial data satisfying $f_{0\eps} \in L^p(\R^d)$ the following basic estimates hold:
\begin{equation}\label{aa17}
\begin{split}
||f_{0\eps}||_{L^{\infty}(0, T; L_+^1 \cap L^{\frac{2d}{d+2}}(\R^d))} &\leq C, \\
||\nabla f_\eps^{(\frac{2d}{d+2}+r-1)/2}||_{L^2(0,T; L^2(\R^d))} &\leq C, \quad 1 < r \leq \frac{2d}{d+2},\\
||f_\eps||_{L^{\frac{3d+2}{d+2}}}(0, T; L^{\frac{3d+2}{d+2}}(\R^d)) &\leq C.
\end{split}
\end{equation}
Applying the weak Young inequality \cite{L-L}
\begin{equation}\label{aa18}
\begin{split}
\int_0^T ||\nabla c_\eps||_{L^2(\R^d)}^{\frac{3d+2}{d+2}} \, dt &\leq C \int_0^T ||(J_\eps \ast f_\eps) \ast \frac{1}{|x|^{d-1}}||_{L^2(\R^d)}^{\frac{3d+2}{d+2}} \, dt \\
&\leq C \int_0^T  ||f_\eps||_{L^{\frac{2d}{d+2}}(\R^d)}^{\frac{3d+2}{d+2}} ||\frac{1}{|x|^{d-1}}||_{L_w^{\frac{d}{d-1}}(\R^d)}^{\frac{3d+2}{d+2}} \, dt \\
&\leq C \int_0^T ||f_\eps||_{L^{\frac{2d}{d+2}}}^{\frac{3d+2}{d+2}} \, dt \leq C,
\end{split}
\end{equation}
So, there exists a subsequence $f_\eps$(without relabeling) such that for any $T>0,$
\begin{equation}\label{aa19}
\begin{split}
f_\eps &\rightharpoonup f \quad \text{in} \quad L^{\frac{3d+2}{d+2}}(0, T; L^{\frac{3d+2}{d+2}}(\R^d)),\\
f_\eps &\overset{\ast}{\rightharpoonup} f \quad \text{in} \quad L^{\infty}(0, T; L_{+}^1 \cap L^{\frac{2d}{d+2}}(\R^d)),\\
\nabla c_\eps &\overset{\ast}{\rightharpoonup} \nabla c  \quad \text{in} \quad L^{\frac{3d+2}{d+2}}(0, T; L^2(\R^d)).\\
\end{split}
\end{equation}
We need to show that the a priori bounds hold uniformly in $\eps$ and we can pass to the limit. We now proceed with the time regularity and applications of the Lions-Aubin Lemma. For any $T > 0:$
\begin{equation}\label{aa20}
\begin{split}
||f_\eps \nabla c_\eps||_{L^2(0, T; L^{\frac{d}{d+1}}(\R^d))} &\leq C, \\
||\nabla f_\eps^{\frac{2d}{d+2}}||_{L^2(0, T; L^{min(2, \frac{4d}{3d+2})}(\R^d))} &\leq C, \\
||\nabla f_\eps ||_{L^{r_2}(0, T; L^{r_2}(\R^d))} &\leq C,
\end{split}
\end{equation}
where $r_2:= min ( 2, \frac{2(\frac{2d}{d+2}+1)}{4-\frac{2d}{d+2}})$. As $f_0 \in L^{\frac{2d}{d+2}}(\R^d),$ by above estimates, one gets for any $T > 0$ and any bounded domain $\Omega,$
\begin{equation}\label{aa21}
||(f_\eps)_t||_{L^{min (2, \frac{2(\frac{2d}{d+2}+1)}{4-\frac{2d}{d+2}})}(0, T; W^{-1, \frac{4d/{d+2}}{\frac{2d}{d+2}+2}}(\Omega))} \leq C.
\end{equation}
Now we proceed with the compactness of $f_\eps$. Let $r_1 = \frac{4d/{d+2}}{\frac{2d}{d+2}+2}$. If $\bar{p}$ satisfies $\frac{dr_1}{d+r_1} \leq \bar{p} < \frac{dr_2}{d-r_2}$, where $r_2$ is defined above, then the following compact embedding holds:
\begin{equation}\label{aa22}
W^{1, r_2}(\Omega) \hookrightarrow \hookrightarrow L^{\bar{p}}(\Omega) \hookrightarrow  W^{-1, r_1}(\Omega)
\end{equation}
this, together with Aubin-Lions Lemma, \eqref{aa20}, \eqref{aa21} implies that

\begin{equation}\label{aa23}
\{f_\eps\}_{\eps>0} \quad \text{is compact in} \quad L^{r_2}(0, T; L^{\bar{p}}(\Omega)).
\end{equation}
Letting $q' = \frac{2(\frac{2d}{d+2}+1)}{4-\frac{2d}{d+2}}$ we see that $\frac{dr_1}{d+r_1} < 2,$  and $\frac{dr_2}{d-r_2} = min (\frac{2d}{d-2}, \frac{dq'}{d-q'}) > 2$ which implies that $\bar{p}=2$ can be chosen. This implies that there exists a subsequence $f_\eps$(not relabeled) such that
\begin{equation}\label{aa24}
f_\eps \to f \quad \text{in} \quad L^{r_2}(0, T; L^{\bar{p}}(\Omega)).
\end{equation}
Let $\{ B_k\}_{k=1}^\infty \subset \R^d$ be a sequence of balls centered at 0 with radius $R_k, \quad R_k \to \infty$. By  a diagonal argument there exists a subsequence (not relabeled) such that the following uniform strong convergence is true:
\begin{equation}\label{aa25}
f_\eps  \to f \quad \text{in} \quad L^{r_2}(0, T; L^{\bar{p}}(B_k)), \forall k \geq 1,
\end{equation}
where $r_2$ and $q'$ are defined above. It now follows from
\begin{equation}\label{aa26}
\begin{split}
||f_\eps||_{L^{\infty}(0, T; L_{+}^1 \cap L^{\frac{2d}{d+2}}(\R^d))} & \leq C,\\
||f_\eps||_{L^{\frac{3d+2}{d+2}}(0, T; L^{\frac{3d+2}{d+2}}(\R^d) )} & \leq C,\\
\end{split}
\end{equation}
and the   H\" older inequality that

\begin{equation}\label{aa27}
\begin{split}
\int_0^T ||f_\eps \nabla c_\eps ||_{\frac{d}{d+2}}^2 \, dt & \leq  \int_0^T ||\nabla c_\eps||_2^2 ||f_\eps||_{\frac{2d}{d+2}}^2 \, dt \\
& \leq C \int_0^T ||f_\eps||_{\frac{2d}{d+2}}^2 \, dt \\
&\leq C(T) \int_0^T ||f_\eps||_{\frac{3d+2}{d+2}}^{{\frac{3d+2}{d+2}}} \, dt \leq C.
\end{split}
\end{equation}
On the other hand, we also have
\begin{equation}\label{aa28}
\begin{split}
||\nabla f_\eps^{\frac{2d}{d+2}}||_{\frac{4d}{3d+2}} & \leq C ||f_\eps^{\frac{1}{2}}||_{\frac{4d}{d+2}}|| \nabla f_\eps^{\frac{3d-2}{2(d+2)}}||_2 \\
&= C  ||f_\eps||_{\frac{2d}{d+2}}^{\frac{1}{2}}||\nabla f_\eps^{\frac{3d-2}{2(d+2)}} ||_2.
\end{split}
\end{equation}
Then, using $||f_\eps||_{L^\infty}(0, T; L^{\frac{2d}{d+2}}) \leq C,$  it follows that $||\nabla f_\eps^{\frac{2d}{d+2}}||_{L^2(0, T; L^{\frac{4d}{3d+2}}(\R^d))} \leq  C.$ Thus,
\begin{equation}\label{aa29}
||\nabla f_\eps^{\frac{2d}{d+2}}||_{L^{2}(0, T; L^{min (2, \frac{4d}{3d+2})})} \leq C.
\end{equation}
By  H\" older inequality
\begin{equation}\label{aa30}
\begin{split}
\int |\nabla f_\eps|^r \, dx  &= C \int |f_\eps^{\frac{3-\frac{4d}{d+2}}{2}}   \nabla f_\eps^{\frac{\frac{4d}{d+2}-1}{2}} |^r \, dx \\
& \leq C    ||f_\eps^{r \frac{3-\frac{4d}{d+2}}{2}}||_{p_1} || |\nabla f_\eps^{\frac{\frac{4d}{d+2}-1}{2}}|^r ||_{q_1},
\end{split}
\end{equation}
where
\begin{equation}\label{aa31}
p_1 = \frac{2\left(\frac{2d}{d+2} +1\right)}{\left[\frac{2(\frac{2d}{d+2} +1)}{4-\frac{2d}{d+2}}\right]\left[3-\frac{4d}{d+2}\right]},  \quad q_1 = \frac{2}{r}, \quad r = \frac{2\left(\frac{2d}{d+2} +1\right)}{4-\frac{2d}{d+2}}.
\end{equation}
Thus, we have for $d < 6,$
\begin{equation}\label{aa32}
\begin{split}
\int_0^T \int |\nabla f_\eps||^r \, dx \, dt  &\leq    C \left(  \int_0^T ||f_\eps^{r \frac{6-d}{2(d+2)}}||_{p_1}^{p_1} \, dt\right)^{1/{p_1}} \left(  \int_0^T || \left| \nabla f_\eps^{\frac{\frac{4d}{d+2}-1}{2}} \right|^r ||_{q_1}^{q_1} \, dt\right)^{1/{q_1}} \\
& \leq C \left(  \int_0^T ||f_\eps||_{\frac{3d+2}{d+2}}^{\frac{3d+2}{d+2}} \right)^{1/{p_1}} \left( \int_0^T || \nabla f_\eps^{\frac{3d-2}{2(d+2)}}   \right)^{1/{q_1}}  \leq C.
\end{split}
\end{equation}
This implies that
\begin{equation}\label{aa33}
||\nabla f_\eps||_{L^{\frac{2(\frac{2d}{d+2}+1)}{4-\frac{4d}{d+2}}}\left(0, T; L^{\frac{2(\frac{2d}{d+2}+1)}{4-\frac{4d}{d+2}}}\left(\R^d\right) \right)} \leq C.
\end{equation}
For $d \geq 3$ taking $r=\frac{d+6}{d+2}$ it follows that
\begin{equation}\label{aa34}
||\nabla f_\eps||_{L^2(0, T; L^2(\R^d))} \leq C.
\end{equation}
Combining the two cases we conclude that
\begin{equation}\label{aa35}
||\nabla f_\eps||_{L^{r_2}(0, T; L^{r_2}(\R^d))} \leq C, \quad r_2 = min \left(2, \frac{2\left(\frac{2d}{d+2}+1\right)}{4-\frac{2d}{d+2}}\right).
\end{equation}
The weak formulation for $f_\eps$ is that $\forall \psi \in C_0^{\infty}(\R^d)$ and any $0<t<\infty,$
\begin{equation}\label{aa36}
\begin{split}
&\int \psi f_\eps(\cdot, t) \, dx - \int \psi f_0 \, dx = \int_0^t \int \Delta \psi \left[  \left(\frac{d-2}{d}\right) ||f_\eps||_{\frac{2d}{d+2}}^{\frac{4}{d+2}}   f_\eps^{\frac{2d}{d+2}} + \eps f_\eps \right] \, dx \, ds  \\
&-\kappa \frac{c_d(d-2)}{2} \int_0^t \iint \frac{\left[ \nabla \psi(x) -\nabla \psi(y)\right] \cdot (x-y)}{|x-y|^2+\eps^2}  \frac{f_\eps(x,s)f_\eps(y,s)}{\left( |x-y|^2 + \eps^2\right)^{{d-2}/2}} \, dx \, dy \, ds.
\end{split}
\end{equation}
In light of the above estimates, we  deduce that  for any bounded domain $\Omega$

\begin{equation}\label{aa37}
f_\eps^{2d/{d+2}}  \longrightarrow  f^{2d/{d+2}} \quad \text{in} \quad L^1(0, T; L^1(\Omega)).
\end{equation}
This directly implies that
\begin{equation}\label{aa38}
\begin{split}
&\int_0^T \int   \Delta \psi \left[  \left(\frac{d-2}{d}\right) ||f_\eps||_{\frac{2d}{d+2}}^{\frac{4}{d+2}} \right]  f_\eps^{\frac{2d}{d+2}} \, dx \, dt \\
& \quad  \longrightarrow \int_0^T \int   \Delta \psi \left[  \left(\frac{d-2}{d}\right) ||f||_{\frac{2d}{d+2}}^{\frac{4}{d+2}}\right]   f^{\frac{2d}{d+2}} \, dx \, dt \quad \text{as} \quad \eps \to 0.
\end{split}
\end{equation}
We  note that  $\frac{1}{|x-y|^d} - \frac{1}{\left( |x-y|^2+\eps^2\right)^{d/2}} \leq \frac{d}{2}\frac{\eps}{|x-y|^{d+1}}$  and  so
\begin{equation}\label{aa39}
\begin{split}
&\qquad \Bigg| \int_0^T \iint \left[  \nabla \psi(x) - \nabla \psi(y)\right]  \cdot (x-y)   \\
& \qquad  \quad \times \left(   \frac{1}{|x-y|^d} - \frac{1}{\left(|x-y|^2+\eps^2\right)}  \right)^{d/2}  f_\eps(x,t) f_\eps(y,t) \, dy \, dx \, dt \Bigg| \\
& \qquad \quad  \leq C\eps \int_0^T \iint   \frac{f_\eps(x,t) f_\eps(y,t)}{|x-y|^{d-1}}  \, dy \, dx \, dt  \\
&  \qquad \qquad   \leq C \eps \int_0^T   ||f_\eps||_{\frac{2d}{d+1}}^2 \, dt \leq C(T) \eps.
\end{split}
\end{equation}
In addition, $\forall \psi \in C_0^{\infty}(\R^d)$,
\begin{equation}\label{aa40}
\begin{split}
&\Bigg|  \iint \left[ \nabla \psi(x) - \nabla \psi(y)\right] \cdot (x-y)    \left(  \frac{f_\eps(x)f_\eps(y)}{|x-y|^d} - \frac{f(x)f(y)}{|x-y|^d}   \right) \, dy \, dx \Bigg|\\
&   \quad \leq C \iint_{\Omega \times \Omega}  \frac{|f_\eps(x)-f(x)|f_\eps(y)}{|x-y|^{d-2}} \, dy \, dx  \\
& \qquad  + C \iint_{\Omega \times \Omega} \frac{|f_\eps(y)-f(y)|f(x)}{|x-y|^{d-2}} \, dy \, dx =: I_1 + I_2.
\end{split}
\end{equation}
For $\int_0^T I_1 \, dt $ we take $\bar{p} = 2$ in the above estimates and deduce that
\begin{equation}\label{aa41}
\begin{split}
&  \int_0^T I_1 \, dt   \\
                            &= C \int_0^T \int_\Omega |f_\eps(x) f(x)|  \left[ \int_{\R^d}  \frac{f_\eps(y)}{|x-y|^{d-2}}  \right] \, dx \, dt \\
		         & \leq  C \int_0^T ||f_\eps-f||_{L^{2d/{d+2}}(\Omega)} ||f_\eps||_{L^{2d/{d+2}}(\R^d)} \, dt \\
 		        &\leq C \left[ \int_0^T  ||f_\eps - f||_{L^2(\Omega)}^{\theta (\frac{3d+2}{2d})}  \right]^{{2d}/{3d+2}} \left[ \int_0^T ||f_\eps||_2^{\theta  (\frac{3d+2}{2d}) } \right]^{{d+2}/{3d+2}},
\end{split}
\end{equation}
where    $\theta = \frac{d-2}{d}$ and  $\theta \left(\frac{3d+2}{2d}\right)  < r_2 := min \left( 2, \frac{2(\frac{2d}{d+2}+1)}{4-\frac{2d}{d+2}} \right)$
and  $\Omega :=  \{  x \in \R^d : f_s(x) > 0\},$ where $f_s(x)$ is the steady solution of \eqref{Meqn}, \eqref{nabc} and \eqref{ic}. $\int_0^T I_2 \, dt$ is estimated in the same manner so we skip the details. Taking  $\eps \to 0 $ in the above estimates we deduce that
\begin{equation}\label{aa42}
\begin{split}
&\int_0^T \iint \left[ \nabla \psi(x) -\nabla \psi(y)  \right]\cdot (x-y)   \\
& \qquad \qquad  \times  \left(  \frac{f_\eps(x)f_\eps(y)}{(|x-y|^2+\eps^2)^{d/2}} - \frac{f(x)f(y)}{|x-y|^d}   \right)  \, dy \, dx \, dt  \longrightarrow  0.
\end{split}
\end{equation}
So, taking the limit $\eps \to 0$ in above estimates we conclude that for any $0 < t < T$,
\begin{equation}\label{aa43}
\begin{split}
& \int_{\R^d}  \psi f(\cdot, t) \, dx - \int_{\R^d} \psi f_(x) \, dx \\
& =  \int_0^t \int_{\R^d} \Delta \psi \left[  \left(\frac{d-2}{d}\right) ||f||_{\frac{2d}{d+2}}^{\frac{4}{d+2}}   f^{\frac{2d}{d+2}}  \right] \, dx \, ds  \\
& - \frac{\kappa c_d }{2} (d-2)  \int_0^t \iint \frac{\left[ \nabla \psi(x) -\nabla \psi(y)\right] \cdot (x-y)}{|x-y|^2+\eps^2}  \frac{f(x,s)f(y,s)}{\left( |x-y|^2 + \eps^2\right)^{{d-2}/2}} \, dx \, dy \, ds.
\end{split}
\end{equation}
It is now possible to make simple adaptations of  the arguments  in Steps 12, 13, 14 and 15 in the proof of the Theorem 2.11 of \cite{BL} to conclude

\begin{itemize}

\item  Strong convergence for the weak solution.

\item  Convergence of the free energy.

\item  Lower semi-continuity of the dissipation term in the energy inequality \eqref{a9}.

\item  Existence of weak entropy solution with   the  energy inequality

\end{itemize}

We then have the following strong convergence:

\begin{equation}\label{aa44}
f_\eps \to f  \quad \text{ in } \quad L^{r_2}(0, T; L^{r_0}(\R^d)),
\end{equation}
where
\begin{equation}\label{aa45}
\begin{split}
1 &\leq r_0 < min \left \{  \frac{dr_2}{d-r_2}, \frac{3d+2}{d+2}  \right \}, \\
r_2 &= min  \left \{  2, \frac{3d+2}{d+4}  \right \}.
\end{split}
\end{equation}
The estimates above yield that(again some further details can be provided)
\begin{equation}\label{aa46}
\begin{split}
||\nabla c_\eps(\cdot, t)||_{L^2(\R^d)}^2  &\longrightarrow  ||\nabla c(\cdot, t)||_{L^2(\R^d)}^2 \quad \text{a.e. in} \quad (0, T),\\
||f_\eps(\cdot, t)||_{L^{2d/{d+2}}(\R^d)}     &\longrightarrow  ||f(\cdot, t)||_{L^{\frac{2d}{d+2}}(\R^d)} \quad \text{a.e. in} \quad (0, T).
\end{split}
\end{equation}
Hence, we obtain,
\begin{equation}\label{aa47}
E_\alpha[f_\eps(\cdot, t)]  \longrightarrow   E_\alpha[f(\cdot, t)] \quad \text{a.e. in} \quad (0, T).
\end{equation}
Lower semi-continuity of the dissipation term follows similarly and we skip the details. Thus, the    existence of a weak entropy solution with   the  energy inequality\eqref{a9} has been established.



\section{Acknowledgments} \label{ack}

The work of E. A. Carlen   is partially supported by U.S. N.S.F. grant DMS 1501007. The authors thank the anonymous referee for the valuable suggestions to improve the presentation of the original manuscript.

\end{document}